\documentclass[11pt]{article}
\usepackage{mathrsfs, cite}
\usepackage{amsmath,amsfonts,amssymb,rotating,amsthm}
\usepackage{psfrag,eepic,color}
\usepackage{framed}

\def\qed{\nopagebreak\hfill{\rule{4pt}{7pt}}}
\def\proof{\noindent {\it{Proof.} \hskip 2pt}}
\newtheorem{theo}{Theorem}[section]

\newtheorem{lem}[theo]{Lemma}

\newtheorem{ex}{Example}[section]
\theoremstyle{remark}

\begin{document}
\begin{center}
{\Large \bf Log-concavity of $P$-recursive sequences}
\end{center}

\begin{center}
{\large Qing-hu Hou and Guojie Li}\\[9pt]
School of Mathematics\\
Tianjin University\\
Tianjin 300072, China\\[6pt]
{\tt
qh\_hou@tju.edu.cn}, {\tt 1017233006@tju.edu.cn}
\end{center}

\vspace{0.3cm} \noindent{\bf Abstract.}

We consider the higher order Tur\'an inequality and higher order log-concavity for sequences $\{a_n\}_{n \ge 0}$ such that
\[
\frac{a_{n-1}a_{n+1}}{a_n^2} = 1 + \sum_{i=1}^m  \frac{r_i(\log n)}{n^{\alpha_i}} + o\left( \frac{1}{n^{\beta}} \right),
\]
where $m$ is a nonnegative integer, $\alpha_i$ are real numbers, $r_i(x)$ are rational functions of $x$ and
\[
0 < \alpha_1 < \alpha_2 < \cdots < \alpha_m < \beta.
\]
We will give a sufficient condition on the higher order Tur\'an inequality and the $\ell$-log-concavity for $n$ sufficiently large.
 Many $P$-recursive sequences fall in this frame. At last,  we will give a method to find the  $N$ such that for any $n>N$,  the higher order Tur\'an inequality holds.

\vskip 15pt


\vskip 15pt


\section{Introduction}
The Tur\'{a}n inequalities and the higher order Tur\'{a}n inequalities arise in the study
of Maclaurin coefficients of real entire functions in the Laguerre-P\'{o}lya class \cite{Sze}.
A sequence $\{a_n\}_{n\geq 0}$ of real numbers is said to satisfy the Tur\'{a}n inequalities
or to be log-concave, if
\[
a_n^2-a_{n-1}a_{n+1}\geq 0, \quad \forall\, n\geq 1.
\]
The Tur\'an inequalities are also called Newton's inequality \cite{Wag77, Nic}. For more results on the log-concavity, we refer to \cite{Des, McN, Sta89}.

A sequence $\{a_n\}_{n \ge 0}$ is said to satisfy the higher order Tur\'an inequalities if for all $n \ge 1$,
\begin{equation}\label{Turan}
  4(a_n^2-a_{n-1}a_{n+1})(a_{n+1}^2-a_na_{n+2})-(a_na_{n+1}-a_{n-1}a_{n+2})^2\geq 0.
\end{equation}
Chen, Jia and Wang\cite{Chen} use the Hardy-Ramanujan-Rademacher formula to prove that, when $n\geq95$, the partition function satisfies the higher order Tur\'{a}n inequality. Dimitrov\cite{Dim} observed that for a real entire function $\Psi(x)$ in the Laguerre-P\'{o}lya class, the Maclaurin coefficients satisfy the higher order
Tur\'{a}n inequalities.  Do\v{s}li\'{c} \cite{Dos} proved some combinatorial sequences are log-balanced. Wang\cite{Wang} proved some combinatorial sequences satisfy the higher order Tur\'{a}n inequalities.  Griffin, Ono, Rolen, and Zagier\cite{Gri} give a straightforward way to check whether a given sequence eventually satisfies the higher Tur\'{a}n inequalities.

Let $\varphi$ be the operator given by
\[
\varphi\{a_n\}_{n\geq 0}=\{a_{n+1}^2-a_{n}a_{n+2}\}_{n\geq 0},
\]
and let $\varphi^k$ be the composition of $\varphi$ $k$ times with itself.  If
\begin{equation}
\varphi\{a_n\}_{n\geq 0},\quad \varphi^2\{a_n\}_{n\geq 0},\quad \ldots, \quad \varphi^r\{a_n\}_{n\geq 0}
\end{equation}
are all non-positive sequences, $\{a_n\}_{n \ge 0}$ is said to satisfy the $\ell$-th iterated Tur\'an inequalities \cite{Cso03} or to be  $\ell$-log-concave \cite{Moll}. Clearly, $1$-log-concavity is just log-concavity.

We focus on the behaviour of a sequence $\{a_n \}_{n \ge 0}$ when $n$ is sufficiently large. Similar to the asymptotic  $\ell$-log-concavity given in \cite{Hou}, we say a sequence satisfies the higher order Tur\'an inequalities asymptotically if \eqref{Turan} holds for $n$ sufficiently large. We aim to give a criterion on the asymptotic higher order Tur\'an inequalities for $P$-recursive sequences.   A $P$-recursive sequence of order $d$ satisfies a recurrence relation of the form
\[
a_n=r_1(n)a_{n-1}+r_2(n)a_{n-1}+\cdots+r_d(n)a_{n-d},
\]
where $r_i(n)$ are rational functions of $n$, see \cite[Section 6.4]{Sta}.

By the asymptotic estimation given by Birkhoff and Trjitzinsky \cite{Birk} and developed by Wimp and Zeilberger \cite{Wim}, we see that (Theorem~\ref{P-un}) many $P$-recursive sequences $\{a_n\}_{n \ge 0}$ satisfy
\begin{equation}\label{un}
u_n = \frac{a_{n-1}a_{n+1}}{a_n^2} = 1 + \sum_{i=1}^m  \frac{r_i(\log n)}{n^{\alpha_i}} + o\left( \frac{1}{n^{\beta}} \right),
\end{equation}
where $m$ is a nonnegative integer, $\alpha_i$ are real numbers, $r_i(x)$ are rational functions of $x$ and
\[
0 < \alpha_1 < \alpha_2 < \cdots < \alpha_m < \beta.
\]

With the asymptotic form \eqref{un}, we are able to give a sufficient condition on the asymptotic higher order Tur\'an inequalities in Section~\ref{sec2}. Then in Section~\ref{sec3}, we extend the result of Hou and Zhang \cite{Hou} on asymptotic  $\ell$-log-concavity. In particular, we apply the criterions on $P$-recursive sequences in Section~\ref{sec4}.
At last, we give a method to find the $N$ such that for any $n>N$ the higher order Tur\'an inequalities hold in section~\ref{sec5}.

\section{The asymptotic higher order Tur\'an inequalities}\label{sec2}
In this section, we will give a sufficient condition on the asymptotic higher order Tur\'an inequalities.

We firstly give an estimation on the difference of rational functions of $\log n$.
\begin{lem}\label{dif-log}
Let $r(x)$ be a rational function of $x$ and $K$ be a positive integer. Then
\[
r(\log (n+1)) - r( \log n)
=  \sum_{i=1}^K \frac{r_i(\log n)}{n^i} + o\left( \frac{1}{n^K} \right),
\]
for some rational functions $r_i(x)$. Similarly,
\[
r(\log (n-1)) - r( \log n)
=  \sum_{i=1}^K \frac{\tilde{r}_i(\log n)}{n^i} + o\left( \frac{1}{n^K} \right),
\]
for some rational functions $\tilde{r}_i(x)$. If $r(x)$ is a polynomial, are the $r_i(x), \tilde{r}_i(x)$'s. Moreover,
\[
r_1(x) = - \tilde{r}_1(x) = r'(x) \quad \mbox{and} \quad
r_2(x)=\tilde{r}_2(x)= \frac{r''(x)-r'(x)}{2}.
\]
\end{lem}
\proof
Let
\begin{equation}\label{delta}
\delta = \log(n+1) - \log n = \sum_{k=1}^K (-1)^{k-1} \frac{1}{k n^k} + o \left( \frac{1}{n^K} \right).
\end{equation}
Notice that $\delta$ behaves like $1/n$ and for any rational function $f(x)$ we have
\[
\lim_{x \to \infty} \frac{f(\log x)}{x}=0.
\]
Therefore,
\begin{align*}
r(\log(n+1)) & = r((\log n) + \delta ) \\
& = \sum_{k=0}^K \frac{r^{(k)}(\log n)}{k!} \delta^k + o \left( \delta^K \right),
\end{align*}
where $r^{(k)}(x)$ is the $k$-th derivative of $r(x)$. Expanding $\delta^k$ by \eqref{delta} and collecting the coefficients of $1/n^k$, we derive that
\[
r(\log(n+1)) = r(\log n) + \sum_{k=1}^K \frac{r_k(\log n)}{n^k} + o\left( \frac{1}{n^K} \right),
\]
where $r_k(x)$ are certain rational functions of $x$. In particular, we have
\[
r_1(x) = r'(x) \quad\text{and}\quad r_2(x) = \frac{r''(x)-r'(x)}2,
\]
completing the proof. \qed

Now we are ready to give the main result.
\begin{theo} \label{th-T}
Let $\{a_n\}_{n \ge 0}$ be a sequence such that \eqref{un} holds.
Assume that $\alpha_m- \alpha_1 \ge 1$. If
\begin{quote}
$0<\alpha_1<2$ and $r_1(x)<0$ for $x$ sufficiently large,
\end{quote}
or
\begin{quote}
$\alpha_1=2$ and $r_1(x)<-1$ for $x$ sufficiently large,
\end{quote}
then $\{a_n\}_{n \ge 0}$ satisfies higher order Tur\'an inequalities asymptotically.
\end{theo}

\proof
Dividing $a_n^2a_{n+1}^2$ on both sides of \eqref{Turan}, we see that the higher order Tur\'an inequality is equivalent to
\begin{equation}\label{T-un}
4(1-u_n)(1-u_{n+1})-(1-u_nu_{n+1})^2 \geq 0.
\end{equation}

Denote
\[
\alpha = \alpha_1,  \quad
u_n = 1 + \frac{\xi(n)}{n^{\alpha}},
\]
and
\[
f(n)  = 4(1-u_n)(1-u_{n+1})-(1-u_nu_{n+1})^2.
\]
We have
\begin{equation}\label{xi-Turan}
f(n)=-\frac{t^2(n)+\xi(n)\xi(n+1) \big( 2\xi(n+1)n^\alpha+2\xi(n)(n+1)^\alpha+\xi(n)\xi(n+1) \big)} {n^{2\alpha}(n+1)^{2\alpha}},
\end{equation}
where
\[
t(n)=\xi(n)(n+1)^\alpha-\xi(n+1)n^\alpha.
\]
Noting that
\[
\xi(n)= \sum_{i=1}^m  \frac{r_i(\log n)}{n^{\alpha_i-\alpha}} + o\left( \frac{1}{n^{\beta-\alpha}} \right),
\]
we have
\begin{align*}
t(n) & =\sum_{i=1}^m \left( \frac{(n+1)^\alpha r_i(\log n)}{n^{\alpha_i-\alpha}} - \frac{n^\alpha r_i(\log(n+1))}{(n+1)^{\alpha_i-\alpha}} \right)
+ o\left( \frac{1}{n^{\beta -2\alpha}} \right).
\end{align*}

From Lemma~\ref{dif-log}, we derive that
\begin{align*}
  & \frac{(n+1)^\alpha r(\log n)}{n^\gamma} - \frac{n^\alpha r(\log (n+1))}{(n+1)^\gamma} \\
 = & n^{\alpha-\gamma} r(\log n) \left( \left(1+\frac{1}{n} \right)^\alpha - \left(1+\frac{1}{n} \right)^{-\gamma} \left( 1 + \frac{r(\log(n+1))-r(\log n)}{r(\log n)} \right) \right) \\
=& n^{\alpha-\gamma} r(\log n)\left( \frac{\alpha+\gamma}{n} - \frac{r'(\log n)}{n r(\log n)} +
o \left( \frac{1}{n} \right) \right) \\
= &  n^{\alpha-\gamma-1} r(\log n) \big( \alpha+\gamma+o(1) \big),
\end{align*}
where the last equality holds since $r'(x)/r(x) \to 0$ when $x \to \infty$ for any rational function $r(x)$. Therefore,
\begin{align*}
t(n)& =\sum_{i=1}^m  n^{2\alpha-\alpha_i-1} r_i(\log n) \big( \alpha_i + o(1) \big) + o\left( \frac{1}{n^{\beta - 2\alpha}} \right) \\
    &= \alpha n^{\alpha-1} r_1(\log n)(1 + o(1)),
\end{align*}
since $\beta>\alpha_m \ge \alpha+1$.
Hence,
\[
t^2(n) = \alpha^2 n^{2 \alpha -2} \left( r_1(\log n) \right)^2 (1+o(1)).
\]
On the other hand, we have
\begin{multline*}
  \xi(n) \xi(n + 1) \big( 2\xi(n + 1)n^\alpha + 2\xi(n)(n + 1)^\alpha + \xi(n) \xi(n + 1) \big) \\
  = 4 \left( r_1(\log n) \right)^3  n^\alpha (1 + o(1)).
\end{multline*}

If $\alpha<2$, we have $2 \alpha -2<\alpha$ and thus
\[
f(n) = -4 \left( \frac{r_1(\log n)}{n^\alpha} \right)^3 ( 1 + o(1)),
\]
which is positive for $n$ sufficiently large if $r_1(x)<0$ for $x$ sufficiently large.

If $\alpha=2$, we have
\[
f(n) = -4 \left( \frac{r_1(\log n)}{n^3} \right)^2 (r_1(\log n)+1+o(1)),
\]
which is positive for $n$ sufficiently large if $r_1(x)<-1$ for $x$ sufficiently large.
\qed

\def\lc{{\rm lc\,}}
\noindent {\it Remark.}
Let $r(x) \in \mathbb{R}[x]$ be a rational function of $x$ with real coefficients. Suppose $r(x)=p(x)/q(x)$, where $p(x),q(x)$ are polynomials in $x$. Then
$r(x)<0$ for $x$ sufficiently large if and only if
\[
\lc p/\lc q<0,
\]
where $\lc p$ and $\lc q$ denotes the leading coefficients of $p$ and $q$, respectively.

\begin{ex}
By Theorem~\ref{th-T}, when $u_n$ are of the following form, the corresponding $\{a_n\}_{n \ge 0}$ satisfies higher order Tur\'an inequalities for $n$ sufficiently large.
\[
1 - \frac{1}{n}, \quad  1 - \frac{1}{n \log n}, \quad 1 - \frac{2}{n^2}, \quad 1 - \frac{\log n}{n^2}.
\]
Noting that for $a_n= 1/n!$, we have
\[
u_n = \frac{n}{n+1} = 1 - \frac{1}{n} + \frac{1}{n^2} + o(n^{-2-\delta}), \quad 0<\delta<1,
\]
and hence $\{ a_n \}_{n \ge 0}$ satisfies higher order Tur\'an inequalities for $n$ sufficiently large. In fact,
\[
4(1-u_n)(1-u_{n+1})-(1-u_n u_{n+1})^2 = \frac{4}{(n+1)(n+2)^2} > 0, \quad \forall\, n \ge 0.
\]
\end{ex}

\begin{ex}
We remark that the condition $\alpha_m \ge \alpha_1+1$ is necessary.

 Let $a_0=a_1=1$ and $a_{n+1}=a_n^2 u_n / a_{n-1}$ for $n \ge 1$ where
\[
u_n = \begin{cases}
1 - \frac{1}{n}, & \mbox{$n$ is even}, \\[5pt]
1 - \frac{1}{n} + \frac{1}{n^{4/3}}, & \mbox{$n$ is odd}.
\end{cases}
\]
We have
\[
a_2=1,\ a_3= \frac{1}{2},\ a_4= \frac{1}{36} (6 + 3^{2/3}),\ a_5= \frac{1}{864} (6 + 3^{2/3})^2, \ldots.
\]
Clearly,
\[
u_n = 1 - \frac{1}{n} + o \left( \frac{1}{n} \right).
\]
However,
\[
4(1-u_n)(1-u_{n+1}) - (1-u_n u_{n+1})^2 = - \frac{1}{n^{8/3}} (1 + o(1)),
\]
and hence $\{a_n\}_{n \ge 0}$ does not satisfy the higher order Tur\'an inequalities when $n$ sufficiently large.
\end{ex}

\begin{ex}
Let $a_n=n \log n$, we see that
\[
u_n = 1 + \frac{-1-\frac{1}{\log n} - \frac{1}{(\log n)^2}}{n^2}  + o\left( \frac{1}{n^{3+\delta}} \right),
\]
for any $0<\delta<1$. Hence $\{a_n\}_{n \ge 1}$ satisfies the  higher order Tur\'{a}n inequalities asymptotically.

\end{ex}

\section{Asymptotic  $\ell$-log-concavity}\label{sec3}

Hou and Zhang \cite{Hou} gave a criterion on the  $\ell$-log-concavity of sequences $\{a_n\}_{n \ge 0}$ such that the $r_i(x)$'s in \eqref{un} are all constants. In this section, we will give a similar criterion for the general case where the $r_i(x)$'s are arbitrary rational functions of $x$.

We first give an estimate on the ratio $a_{n-1}a_{n+1}/a_n^2$ for a special kind of $a_n$.
\begin{lem} \label{le-3.2}
Let
\[
a_n = \frac{r(\log n)}{n^\alpha},
\]
where $r(x)$ is a rational function of $x$ and $\alpha$ be a positive real number. Then for any integer $K \ge 2$,
there exist rational functions $r_i(x)$ of $x$ such that
\[
u_n = \frac{a_{n-1}a_{n+1}} {a_n^2} = 1 + \sum_{i=2}^K \frac{r_i(\log n)}{n^i} + o\left( \frac{1}{n^K} \right).
\]
At the same time, there exist rational functions $\tilde{r}_i(x)$ of $x$ such that
\[
a_{n+1}+a_{n-1}-2 a_n = \sum_{i=2}^K \frac{\tilde{r}_i(\log n)}{n^{i+\alpha}} + o\left( \frac{1}{n^{K+\alpha}} \right).
\]
Moreover, we have
\[
r_2(x)= \alpha + (\log r(x))'' - (\log r(x))'
\]
and
\[
\tilde{r}_2(x) = \alpha(\alpha+1) r(x) - (2 \alpha+1) r'(x) +r''(x).
\]

\end{lem}
\proof
Denote $\log n$, $\log (n+1)$, and $\log (n-1)$ by $l$, $l^+$ and $l^-$, respectively. We have
\begin{align*}
u_n & = \frac{n^{2\alpha}}{(n-1)^{\alpha}(n+1)^{\alpha}}\frac{r(l^+) r(l^-)}{r^2(l)} \\
& = \left(1 - \frac{1}{n^2} \right)^{-\alpha} \left(1 + \frac{r(l^+) - r(l)}{r(l)} \right) \left(1 + \frac{r(l^-) - r(l)}{r(l)} \right).
\end{align*}
By Lemma~\ref{dif-log}, we derive that
\begin{align*}
u_n & = \left(1 + \sum_{i=1}^K {-\alpha \choose i} \left(- \frac{1}{n^2} \right)^i + o \left( \frac{1}{n^{K}} \right) \right) \left(1 +  \frac{1}{r(l)} \sum_{i=1}^K \frac{p_i(l)}{n^i} + o \left( \frac{1}{n^{K}} \right) \right)\\
& \qquad  \left(1 +  \frac{1}{r(l)} \sum_{i=1}^K \frac{q_i(l)}{n^i} + o \left( \frac{1}{n^{K}} \right) \right) \\
& = 1 + \frac{\alpha}{n^2} + \frac{r(l)(r''(l)-r'(l))-(r'(l))^2}{n^2 r^2(l)} + \cdots + o \left( \frac{1}{n^K} \right),
\end{align*}
where $p_i(x),q_i(x)$ are rational functions of $x$.

At the same time, we have
\begin{align*}
& a_{n+1}+a_{n-1}-2 a_n \\
& = \frac{r(l)}{n^\alpha} \left( \left(1+\frac{1}{n}\right)^{-\alpha} \left(1+\frac{r(l^+)-r(l)}{r(l)} \right) \right. \\
& \qquad \qquad + \left. \left(1-\frac{1}{n}\right)^{-\alpha} \left(1+\frac{r(l^-)-r(l)}{r(l)} \right) - 2 \right) \\
& = \frac{r(l)}{n^\alpha} \left( \sum_{i=0}^{K} \frac{{-\alpha \choose i}}{n^i} \left(1 +  \frac{1}{r(l)} \sum_{i=1}^{K} \frac{p_i(l)}{n^i} \right) \right. \\
 & \qquad \qquad + \left. \sum_{i=0}^{K} \frac{{-\alpha \choose i}}{(-n)^i} \left(1 +  \frac{1}{r(l)} \sum_{i=1}^{K} \frac{q_i(l)}{n^i} \right)  - 2 + o\left( \frac{1}{n^{K}} \right)\right) \\
& = \frac{r(l)}{n^\alpha} \left( \frac{\alpha(\alpha+1)}{n^2} - \frac{2 \alpha r'(l)}{n^2 r(l)}
+ \frac{r''(l)-r'(l)}{n^2 r(l)} + \cdots + o\left( \frac{1}{n^{K}} \right) \right),
\end{align*}
completing the proof. \qed

Now we are ready to give a criterion on the asymptotic $\ell$-log-concavity.
\begin{theo}\label{th-3.3}
Let $\{a_n\}_{n \ge 0}$ be a sequences such that \eqref{un} holds and denote  $\ell = \lfloor\alpha_m/\alpha_1\rfloor$.
If
\begin{quote}
$0<\alpha_1<2$ and $r_1(x)<0$ for $x$ sufficiently large,
\end{quote}
or
\begin{quote}
$\alpha_1=2$ and $r_1(x)<-2+1/2^{r-2}$ for $x$ sufficiently large,
\end{quote}
then $\{a_n\}_{n \ge 0}$ is asymptotically  $\ell$-log-concave.
\end{theo}
\proof
Let
\[
u_n = \frac{a_{n-1}a_{n+1}}{a_n^2}, \quad
b_n = a_{n}^2 - a_{n-1} a_{n+1}.
\]
It is easy to check that
\[
\frac{b_{n-1}b_{n+1}}{b_n^2} =
u_n^2\frac{(u_{n-1}-1)(u_{n+1}-1)}{(u_n-1)^2}.
\]
We will give an estimation for this ratio.

Denote $\alpha=\alpha_1$ and $r(x)=r_1(x)$.
Since $u_n$ is of form \eqref{un}, so is $u_n^2$. Moreover,
\begin{equation} \label{un2}
u_n^2=1 + \frac{2 r(\log n)}{n^{\alpha}} + \cdots + o\left( \frac{1}{n^{\beta}} \right).
\end{equation}
Rewrite $u_n$ as $u_n = 1 + f_n g_n$ with
\[
f_n=\frac{r(x)}{n^{\alpha}} \quad \mbox{and} \quad
g_n=1 + \sum_{i=2}^{m}\frac{r_i(x)}{r(x)n^{\alpha_i-\alpha}} + o\left(\frac{1}{n^{\beta-\alpha}}  \right).
\]
Then we have
\[
\frac{(u_{n-1}-1)(u_{n+1}-1)}{(u_n-1)^2}=\frac{f_{n-1}f_{n+1}}{f_n^2}\frac{g_{n-1}g_{n+1}}{g_n^2}.
\]
By Lemma~\ref{le-3.2}, we get that
\begin{equation}\label{fn}
\frac{f_{n-1}f_{n+1}}{f_n^2} = 1 +
\frac{t(\log n)}{n^2}+\cdots+ o\left(\frac{1}{n^{\beta}} \right),
\end{equation}
with
\[
t(x)=\alpha + (\log r(x))'' - (\log r(x))'.
\]

To estimate the ratio $g_{n-1}g_{n+1}/g_n^2$, we consider $h_n = \log g_n$. By the Taylor expansion of $\log(1+x)$, we derive that $h_n$ is of form
\[
\sum_{i=1}^l \frac{s_i(\log n)}{n^{\beta_i}} + o \left( \frac{1}{n^{\beta-\alpha}} \right),
\]
where
\[
\alpha_2-\alpha = \beta_1 < \beta_2 < \ldots < \beta_l < \beta-\alpha,
\]
and $s_i(x)$ are rational functions of $x$. Moreover, $s_1(x)=r_2(x)/r(x)$. By Lemma~\ref{le-3.2}, we see that
\begin{align*}
\log\frac{g_{n+1}g_{n-1}}{g_n^2}&=h_{n+1} + h_{n-1} - 2 h_n \\
 &=\sum_{i=1}^l \sum_{j=2}^K \frac{s_{ij}(\log n)}{n^{j+\beta_i}} + o \left( \frac{1}{n^{\beta-\alpha}} \right),
\end{align*}
where $K$ is an integer with $K>\beta-\alpha$ and $s_{ij}(x)$ are rational functions of $x$. By the Taylor expansion of $e^x$, we see that
\begin{equation}\label{gn}
\frac{g_{n+1}g_{n-1}}{g_n^2} = 1 + \sum_{i=1}^{d} \frac{\tilde{s}_i(\log n)}{n^{\gamma_i}} + o \left( \frac{1}{n^{\beta-\alpha}} \right),
\end{equation}
where
\[
\alpha_2-\alpha + 2 = \gamma_1 < \gamma_2 < \ldots < \gamma_{d} < \beta-\alpha,
\]
and $\tilde{s}_i(x)$ are rational functions of $x$.

Combining \eqref{un2}, \eqref{fn} and \eqref{gn} together, we derive that
$b_{n+1}b_{n-1}/b_n^2$ is of form \eqref{un} and
\begin{multline}\label{guji}
\frac{b_{n+1}b_{n-1}}{b_n^2}
= \left( 1 +\frac{2r( \log n)}{n^\alpha}+\cdots+o \left( \frac{1}{n^{\beta}} \right) \right) \\
\times \left( 1 + \frac{t(\log n)}{n^2} + \cdots + o \left( \frac{1}{n^{\beta-\alpha}} \right) \right).
\end{multline}

When $\alpha<2$, we have
\[
\frac{b_{n+1}b_{n-1}}{b_n^2} = 1 + \frac{2r( \log n)}{n^\alpha} + \cdots + o \left( \frac{1}{n^{\beta-\alpha}} \right).
\]
If $r(x)<0$ for $x$ sufficiently large, so is $2r(x)$. If $\alpha_m \ge 2\alpha$, we have $\beta-\alpha > \alpha$. Hence the ratio $b_{n-1}b_{n+1}/b_n^2$ is less than $1$ for $n$ sufficiently large, implying that $\{a_n\}_{n \ge 0}$ is asymptotically $2$-log-concave. Repeating this argument, we finally derive that
$\{a_n\}_{n \ge 0}$ is  $\ell$-log-concave.

Then we consider the case of $\alpha=2$. It is clear that if $r(x)<0$ for $x$ sufficiently large, then $\{a_n\}_{n \ge 0}$ is log-concave. Now suppose that $r(x)<-1$ for $x$ sufficiently large and $\alpha_m \ge 4$. We have
\[
\frac{b_{n+1}b_{n-1}}{b_n^2}
= 1 + \frac{2 r( \log n) + t(\log n)}{n^2} + \cdots + o \left( \frac{1}{n^{\beta-2}} \right).
\]
If
\[
\lim_{x \to \infty} r(x) = a < -1,
\]
we have
\[
\lim_{x \to \infty} 2r(x)+t(x) = 2a+2<0,
\]
since
\[
\lim_{x \to \infty} (\log r(x))' = \lim_{x \to \infty} \frac{r'(x)}{r(x)} = 0,
\]
and
\[
\lim_{x \to \infty} (\log r(x))'' = \lim_{x \to \infty} \frac{(r'(x)/r(x))'}{r'(x)/r(x)} \lim_{x \to \infty} \frac{r'(x)}{r(x)} = 0.
\]
If
\[
\lim_{x \to \infty} r(x) = -1,
\]
we have
\begin{align*}
\lim_{x \to \infty} \frac{(\log r(x))''-(\log r(x))'}{2 r(x) + 2}  & =
\lim_{x \to \infty} \frac{r'(x)/r(x)}{2 r(x) + 2} \lim_{x \to \infty}\frac{(\log r(x))''-(\log r(x))'}{r'(x)/r(x)} \\
& =
\lim_{x \to \infty} \frac{(r'(x)/r(x))'}{2 r'(x)} \cdot (- 1) = 0.
\end{align*}
Hence
\[
2 r(x)+t(x) = (2r(x)+2)(1 + o(1)) < 0
\]
for $x$ sufficiently large, implying that $\{a_n\}_{n \ge 0}$ is $2$-log-concave for $n$ sufficiently large.

In general, if $r(x)<-2+1/2^{\ell-2}$ for $x$ sufficiently large, we have
\[
2 r(x)+t(x) < -2 + \frac{1}{2^{\ell-3}}
\]
for $n$ sufficiently large. Repeating this discussion, we finally derive that $\{a_n\}_{n \ge 0}$ is asymptotically  $\ell$-log-concave. \qed

We see that the condition for the higher order Tur\'an inequalities given by Theorem~\ref{th-T} and the condition for the $2$-log-concavity given by Theorem~\ref{th-3.3} coincide.

\begin{ex}
Let $a_n=n^\alpha (\alpha \ge 2)$, we see that
\[
u_n = \left(1 - \frac{1}{n^2} \right)^\alpha
=1-\frac{\alpha}{n^2}+\sum_{i=2}^{\infty}\frac{{\alpha \choose i}(-1)^i}{n^{2i}}.
\]
Hence $\{a_n\}_{n \ge 0}$ is asymptotically  $\ell$-log-concave for any positive integer  $\ell$.
\end{ex}

\begin{ex}
Let $a_n=n^2 \log n$, we see that
\[
u_n =1 + \frac{-2-\frac{1}{\log n} - \frac{1}{(\log n)^2}}{n^2}  + \sum_{i=2}^\infty \frac{r_i(\log n)}{n^{2i}},
\]
where $r_i(x)$ are rational functions of $x$.
Hence $\{a_n\}_{n \ge 1}$ is asymptotically  $\ell$-log-concave for any positive integer  $\ell$.
\end{ex}

\section{$P$-recursive sequences and asymptotic higher order Tur\'an property}
\label{sec4}

In this section, we apply the criterion given in the previous section to $P$-recursive sequences. Recall that a $P$-recursive sequence of order $d$ satisfies a recurrence relation of the form
\[
a_n=r_1(n)a_{n-1}+r_2(n)a_{n-1}+\cdots+r_d(n)a_{n-d},
\]
where $r_i(n)$ are rational functions of $n$, see \cite[Section 6.4]{Sta}.
We will use the asymptotic estimation of $P$-recursive sequences given by Birkhoff and Trjitzinsky \cite{Birk} and developed by Wimp and Zeilberger \cite{Wim}. They showed that a $P$-recursive sequence is asymptotically equal to a linear combination of terms of the form
\begin{equation} \label{Qs}
e^{Q(\rho,n)} s(\rho,n),
\end{equation}
where
\[
Q(\rho,n)=\mu_0 n \log n + \sum_{j=1}^{\rho}\mu_jn^{j/\rho},
\]
and
\[
s(\rho,n)=n^r\sum_{j=0}^{t-1}(\log n)^j\sum_{s=0}^{M-1}b_{sj}n^{-s/\rho},
\]
with $\rho$, $t$, $M$ being positive integers and $\mu_j$, $r$, $b_{sj}$ being complex numbers.

The following theorem shows that the $u_n$ of a $P$-recursive sequence is of form \eqref{un}.
\begin{theo}\label{P-un}
Let $\{a_n\}_{n \ge 0}$ be a $P$-recursive sequence with asymptotic form \eqref{Qs}.
Then $u_n = a_{n-1}a_{n+1}/a_n^2$ is of the form \eqref{un} when $n$ tends to infinity.
\end{theo}
\proof
Noting that the product of two terms of the form \eqref{un} is still of the form \eqref{un}, we can treat each factor of \eqref{Qs} separately. The factors $Q(\rho,n)$ and $n^r$ has been discussed in \cite[Theorem 3.1]{Hou}. We need only consider the term
\[
f_n = \sum_{s=0}^{M} \sum_{j=0}^{t-1} b_{sj} \frac{(\log n)^j}{n^{s/\rho}} + o \left( \frac{1}{n^{\beta}} \right),
\]
where $b_{0j}$ are not all zeros and $\beta$ is a real number satisfying
\[
(M+1)/\rho>\beta>M/\rho.
\]

Rewrite $f_n$ as
\[
f_n = \sum_{s=0}^M \frac{p_s(\log n)}{n^{s/\rho}} + o \left( \frac{1}{n^{\beta}} \right),
\]
where
\[
p_s(x) = \sum_{j=0}^{t-1} b_{sj} x^j
\]
are polynomials in $x$. By Lemma~\ref{dif-log}, we have
\begin{align*}
f_{n+1} &= \sum_{s=0}^M \frac{p_s(\log (n+1))}{(n+1)^{s/\rho}} + o \left( \frac{1}{n^{\beta}} \right)\\
&= \sum_{s=0}^M \frac{p_s(\log n) + \sum_{j=1}^M \frac{q_{sj}(\log n)}{n^j} + o\left( \frac{1}{n^M} \right)}{n^{s/\rho}}  \left( 1+ \frac{1}{n} \right)^{-s/\rho}  + o \left( \frac{1}{n^{\beta}} \right) \\
&= f_n + \sum_{s=0}^M \frac{p'_s(\log n) - s p_s(\log n)/\rho}{n^{1+s/\rho}} + \sum_{s=0}^M \frac{r_{s}(\log n)}{n^{2+s/\rho}}  + o \left( \frac{1}{n^{\beta}} \right),
\end{align*}
where $q_{sj}(x), r_s(x)$ are polynomials in $x$. Similarly, we have
\begin{align*}
f_{n-1} &= \sum_{s=0}^M \frac{p_s(\log (n-1))}{(n-1)^{s/\rho}} + o \left( \frac{1}{n^{\beta}} \right)\\
&= \sum_{s=0}^M \frac{p_s(\log n) + \sum_{j=1}^M \frac{\tilde{q}_{sj}(\log n)}{n^j} + o\left( \frac{1}{n^M} \right)}{n^{s/\rho}}  \left( 1 - \frac{1}{n} \right)^{-s/\rho}  + o \left( \frac{1}{n^{\beta}} \right) \\
&= f_n - \sum_{s=0}^M \frac{p'_s(\log n) - s p_s(\log n)/\rho}{n^{1+s/\rho}} + \sum_{s=0}^M \frac{\tilde{r}_{s}(\log n)}{n^{2+s/\rho}} + o \left( \frac{1}{n^{\beta}} \right),
\end{align*}
where $\tilde{q}_{sj}(x), \tilde{r}_s(x)$ are polynomials in $x$.
Since $p_0(x) \not = 0$, we have
\[
\frac{1}{f_n^2} \cdot o \left( \frac{1}{n^{\beta}} \right) = o \left( \frac{1}{n^{\beta'}} \right)
\]
for any $\beta'<\beta$. Therefore,
\begin{align*}
\frac{f_{n+1}f_{n-1}}{f_n^2}
& = 1 - \frac{1}{f_n^2}
\left(\sum_{s=0}^M \frac{p'_s(\log n) - s p_s(\log n)/\rho}{n^{1+s/\rho}} \right)^2
+ \frac{1}{f_n^2} \sum_{s=0}^M \frac{t_{s}(\log n)}{n^{2+s/\rho}} \\
& \qquad + o \left( \frac{1}{n^{\beta'}} \right) \\
& = 1 +  \sum_{s=0}^M \frac{\tilde{t}_{s}(\log n)}{n^{2+s/\rho}} + o \left( \frac{1}{n^{\beta'}} \right)
\end{align*}
where $t_s(x)$ are polynomials in $x$ and $\tilde{t}_s(x)$ are rational functions of $x$. \qed

We give some examples.
\begin{ex}
Let
\[
C_n = \frac{1}{n+1} {2n \choose n}
\]
be the Catalan numbers. Its inverse $a_n=1/C_n$ satisfies
\[
 (4n+2)a_{n+1} - (n+2)a_n =0.
\]
By packages such as {\tt Asyrec} given by Zeilberger \cite{Zeil}, {\tt asymptotics.m} given by Kauers \cite{Kau}, or {\tt P-rec.m} given by Hou and Zhang \cite{Hou}, we find that
\[
a_n=c \cdot \frac{n^{3/2}}{4^n} \left( 1 + \frac{9}{8n} + \frac{17}{128 n^2}  +  o\left( \frac{1}{n^2} \right) \right),
\]
where $c$ is a constant. (In fact, $c=\sqrt{\pi}$.)
Therefore,
\[
u_n=1 - \frac{3}{2 n^2} + \frac{9}{ 4n^3} - \frac{21}{8n^4} + o\left( \frac{1}{n^4} \right).
\]
Hence $\{1/C_n\}_{n \ge 0}$ satisfies the  higher order Tur\'{a}n inequalities asymptotically and is asymptotically $2$-log-concave.
\end{ex}

\begin{ex}
Let $I_n$ be the number of involutions on $\{1, \ldots, n\}$, i.e., permutations $\sigma$ such that $\sigma^2={\rm id}$. It is well-known that
\[
I_n = I_{n-1}+(n-1)I_{n-2}.
\]
Hence, $a_n = I_n/n!$ satisfies the recurrence relation
\[
n a_n = a_{n-1}+a_{n-2}.
\]
By the above packages, we find that
\[
a_n = c \cdot \frac{e^{n/2+\sqrt{n}}}{n^{n/2}} \frac{1}{\sqrt{n}} \left( 1 + \frac{7}{24 \sqrt{n}} - \frac{215}{1152 n} + o\left( \frac{1}{n} \right) \right).
\]
Therefore,
\[
u_n = 1 - \frac{1}{2 n} - \frac{1}{4 n^{3/2}} + \frac{5}{8 n^2} +  o\left( \frac{1}{n^2} \right),
\]
and thus $\{I_n/n!\}_{n \ge 0}$ satisfies the  higher order Tur\'{a}n inequalities asymptotically and is asymptotically  $\ell$-log-concave for any positive integer  $\ell$.

Let $A_n$ be the $n$-th Ap\'ery number given by
\[
A_n = \sum_{k=0}^n {n \choose k}^2 {n+k \choose k}^2.
\]
By a similar discussion, we derive that $A_n/n!$ satisfies the higher order Tur\'{a}n inequalities  asymptotically and is asymptotically  $\ell$-log-concave for any positive integer  $\ell$.
\end{ex}

\section{The  higher order Tur\'an inequality}\label{sec5}
In this section, we will give a method to find the explicit $N$ such that the higher order Tur\'an inequalities holds for
$\{a_n\}_{n\geq N}$ when $\{a_n\}_{n\geq0}$ is $P$-recursive.

The following lemma indicates that to prove the higher order Tur\'an inequality, we need only consider the boundary points.
\begin{lem}\label{lem4}
Let
\[
 t(x,y)=4(1-x)(1-y)-(1-xy)^2.
\]
Suppose that there exist $x_1<x_2$ and $y_1<y_2$ such that
\[
t(x_1,y_1)>0, \quad t(x_1,y_2)>0, \quad t(x_2,y_1)>0, \quad t(x_2,y_2)>0.
\]
Then for every $x$ and $y$ such that
\[
x_1<x<x_2, \quad  y_1<y<y_2,
\]
we have $t(x,y)>0$.
\end{lem}
\proof
It's easy to see
 \[
 t(x,y_1)=-y_1^2x^2+(6y_1-4)x+3-4y_1.
 \]
 Since $-y_1^2\leq0$ and $t(x_1,y_1)>0$, $t(x_2,y_1)>0$,   for any $x\in (x_1,x_2)$, we have $t(x,y_1)>0$.
 Similarly, $t(x,y_2)>0$ for any $x\in (x_1,x_2)$.

 We also can rewrite $t(x,y)$ as
 \[
 t(x,y)=-x^2y^2+(6x-4)y+3-4x.
 \]
 Since $-x^2\leq0$ and  $t(x,y_1)>0$, $t(x,y_2)>0$, we have $t(x,y)>0$ for any $y\in (y_1,y_2)$. \qed

 As a natural consequence, we derive that

\begin{theo}\label{th5.2}
Let
\[
t(x,y)=4(1-x)(1-y)-(1-xy)^2 \quad \mbox{and} \quad
u_n = \frac{a_{n-1}a_{n+1}}{a_n^2}.
\]
If there exist an integer $N$, an upper bound $f_n$ and a lower bound $g_n$ of $u_n$ such that for all $n \ge N$,
\[
g_n<u_n<f_n,
\]
and
\[
t(g_n,g_{n+1})>0, \quad t(g_n,f_{n+1})>0, \quad
t(f_n,g_{n+1})>0, \quad t(f_n,f_{n+1})>0.
\]
Then  $\{a_n\}_{n\geq N}$ satisfies the higher order Tur\'an inequality.
\end{theo}

When $\{a_n\}$ has the asymptotic expression ~(\ref{Qs}) such that $t=1$, we are able to compute the bounds $f_n$ and $g_n$ by the following HT algorithm. We implemented the HT algorithm in the package {\tt P}-rec.m
which is available at \cite{HT}.
Here we give an outline of the  algorithm.
\begin{framed}
\begin{center}
{\Large \bf The HT Algorithm}
\end{center}

INPUT: a difference operator  $L$ which annihilates $a_n$ ,  initial values $a_0,a_1,\ldots,a_m$, an integer $K$.

OUTPUT: An integer $N$, rational functions $g_n$ and $f_n$ such that
\[
g_n<\frac{a_{n-1}a_{n+1}}{a_n^2}<f_n, \quad n>N.
\]

1.  Find the asymptotic expansion of $a_n$ up to $K$ terms.

2. Denote $r_n=\frac{a_{n}}{a_{n-1}}$, compute the asymptotic expression of   $r_n$
\[
r_n \approx \sum_{i=1}^{K}\frac{c_i}{n^{\alpha_i}},
\]
where $c_i$, $\alpha_i$ are real numbers .

3. Use the  method in \cite{Hou} to find $N_1$ such that
\[
s_l(n)<\frac{a_{n}}{a_{n-1}}<s_u(n), \quad n>N_1.
\]
where
\[
s_u(n)=\sum_{i=1}^{K}\frac{c_i}{n^{\alpha_i}}+\frac{1}{n^{\alpha_K}}, \quad  s_l(n)=\sum_{i=1}^{K}\frac{c_i}{n^{\alpha_i}}-\frac{1}{n^{\alpha_K}}.
\]

4. Denote $u_n=\frac{a_{n+1}a_{n-1}}{a_n^2}$, compute the asymptotic expression of $u_n$
\[
u_n\thickapprox 1+\sum_{i=1}^{K-1}\frac{d_i}{n^{\beta_i}},
\]
where $d_i$ are real numbers and $\beta_i$ are positive real numbers.

5. Let
\[
f_n=1+\sum_{i=1}^{K-1}\frac{d_i}{n^{\beta_i}}+\frac{1}{n^{\beta_{K-1}}}, \quad g_n=1+\sum_{i=1}^{K-1}\frac{d_i}{n^{\beta_i}}-\frac{1}{n^{\beta_{K-1}}}.
\]
Find $N_2$ such that
\[
\frac{s_u(n+1)}{s_l(n)}<f_n, \quad \frac{s_l(n+1)}{s_u(n)}>g_n, \quad n>N_2.
\]
6. Return
\[
  \{g_n, \quad f_n, \quad \max\{N_1,N_2\}\}.
\]
\end{framed}

Now I'd like to show the usage of the package.
\begin{ex}
Suppose $\{a_n\}$ is defined by
\[
a_n=\sum_{k=0}^{n}\binom{n}{k}^4.
\]
Then by Zeilberger's algorithm we  find that $\{a_n\}$ satisfies the  recurrence relation:
\[
(n+2)^3a_{n+2}-2(3n^2+9n+7)(2n+3)a_{n+1}-4(n+1)(4n+3)(4n+5)a_n=0,
\]
with the initial values
\[
a_0=1, \quad a_1=2, \quad a_2=18.
\]
By the command
\[
HT[(n+2)^3N^2-2(3n^2+9n+7)(2n+3)N-4(n+1)(4n+3)(4n+5),n,N,\{1,2,18\},4]
\]
We  get
\[
\{\{1+\frac{1}{2n^2},1+\frac{5}{2n^2}\},94\},
\]
which means
\[
1+\frac{1}{2n^2}<\frac{a_{n+1}a_{n-1}}{a_n^2}<1+\frac{1}{2n^2}, \quad n\geq94.
\]
\end{ex}

Now we use  Theorem ~\ref{th5.2}  to prove the higher order Tur\'an inequalities given in \cite{Wang}.
 Our method can treat those $P$-recursive sequences of order $2$ directly. Moreover, we can treat $P$-recursive sequences of higher order.
\begin{ex}
(Theorem 3.2 of \cite{Wang}) Let $M_n$ be the Motzkin numbers given by
\[
(n+4)M_{n+2}-(2n+5)M_{n+1}-3(n+1)M_n=0, \quad n\geq0,
\]
with the initial values
\[
M_0=1, \quad M_1=1.
\]
Then the sequence $\{\frac{M_n}{n!}\}_{n\geq1}$ satisfies the higher order Tur\'an inequality.
\end{ex}
\proof
By  the algorithm HT, we derive
\[
1+\frac{1}{2n^2}<\frac{M_{n+1}M_{n-1}}{M_n^2}<1+\frac{5}{2n^2}, \quad n\geq75.
\]
Let
\[
a_n=\frac{M_n}{n!}, \quad  u_n=\frac{a_{n+1}a_{n-1}}{a_n^2}.
\]
Then we have
\[
\frac{n}{n+1} \left( 1+\frac{1}{2n^2} \right) <u_n<\frac{n}{n+1} \left( 1+\frac{5}{2n^2} \right) , \quad n\geq75.
\]
Let
\[
g_n=\frac{n}{n+1} \left( 1+\frac{1}{2n^2} \right) \quad \mbox{and} \quad
f_n=\frac{n}{n+1} \left( 1+\frac{5}{2n^2} \right).
\]
Direct computation shows that
\[
t(g_n,g_{n+1})=\frac{64n^5+96n^4+16n^3-24n^2-8n-9}{16n^2(n+1)^4(n+2)^2},
\]
which is positive for $n>0$.

Similarly, we have
\[
t(g_n,f_{n+1})=\frac{64n^5-224n^4-240n^3+24n^2+152n-49}{16n^2(n+1)^4(n+2)^2},
\]
which is positive for $n>4$.
\[
t(f_n,g_{n+1})=\frac{64n^5-96n^4-496n^3-680n^2-520n-225}{16n^2(n+1)^4(n+2)^2},
\]
which is positive for $n>4$.
\[
t(f_n,f_{n+1})=\frac{64n^5-288n^4+272n^3+264n^2-360n-1225}{16n^2(n+1)^4(n+2)^2},
\]
which is positive for $n>3$.

Checking the initial values,  we finally get that $\{\frac{M_n}{n!}\}_{n \ge 1}$ satisfies the higher order Tur\'an inequality. \qed

\begin{ex}
(Theorem 3.12 of \cite{Wang}) Let $F_n^{(3)}$ be the Franel numbers of order $3$ given by
\[
(n+2)^2F_{n+2}^{(3)}-(7n^2+21n+16)F_{n+1}^{(3)}-8(n+1)^2F_n^{(3)}=0, \quad n\geq0,
\]
with the initial values
\[
F_0^{(3)}=1, \quad F_1^{(3)}=2.
\]
 Then the sequence $\{\frac{F_n^{(3)}}{n!}\}_{n\geq1}$ satisfies the higher order Tur\'an inequality.
\end{ex}
\proof
By the algorithm HT, we derive
\[
1<\frac{F_{n+1}^{(3)}F_{n-1}^{(3)}}{(F_n^{(3)})^2}<1+\frac{2}{n^2}, \quad n\geq16.
\]
Let
\[
a_n=\frac{F_n^{(3)}}{n!}, \quad  u_n=\frac{a_{n+1}a_{n-1}}{a_n^2}.
\]
Then we have
\[
\frac{n}{n+1}<u_n<\frac{n}{n+1}\left(1+\frac{2}{n^2}\right), \quad n\geq16.
\]
Let
\[
g_n=\frac{n}{n+1}, \quad \mbox{and} \quad
f_n=\frac{n}{n+1}\left(1+\frac{2}{n^2}\right).
\]
Direct computation shows that
\[
t(g_n,g_{n+1})=\frac{4}{(n+1)(n+2)^2},
\]
which is positive for $n\geq0$.

Similarly, we have
\[
t(g_n,f_{n+1})=\frac{4n^3-8n^2-20n-12}{(n+1)^4(n+2)^2}>0, \quad n>3.
\]

\[
t(f_n,g_{n+1})=\frac{4n^2-12n-4}{n^2(n+1)(n+2)^2}>0,  \quad n>3.
\]

\[
t(f_n,f_{n+1})=\frac{4n^5-12n^4+4n^3+12n^2-8n-36}{n^2(n+1)^4(n+2)^2}>0, \quad n>2.
\]

Checking the initial values,  we finally get that $\{\frac{F_n^{(3)}}{n!}\}_{n\geq1}$ satisfies the higher order Tur\'an inequality. \qed

Use the same method, we also can prove that $\{\frac{F_n}{n!}\}_{n\geq3}$ and $\{\frac{D_n}{n!}\}_{n\geq0}$ satisfy the higher order Tur\'an inequality, where $F_n$ is
the Fine number and $D_n$ is the Domb number.  These properties have been proved in \cite{Wang}.

\begin{ex}
Consider the sequences $\{b_n\}_{n\geq0}$ given by
\[
b_n=\sum_{k=0}^{n}\frac{\binom{n}{k}\binom{n+k}{k}}{2k-1}.
\]
Then the sequence $\{\frac{b_n}{n!}\}_{n\geq0}$ satisfies the higher order Tur\'an inequality.
\end{ex}
\proof
Use Zeilberger's algorithm we see that $\{b_n\}_{n\geq0}$ satisfies the following recurrence relation
\[
(n+3)b_{n+3}-(7n+13)b_{n+2}+(7n+15)b_{n+1}-(n+1)b_n=0, \quad n\geq0,
\]
with the initial values
\[
b_0=-1,\quad b_1=1,\quad b_2=7.
\]
By the algorithm HT, we have the following estimates:
\[
1+\left(\frac{3(17+12\sqrt{2})}{2(3+2\sqrt{2})^2} - 1 \right) \frac{1}{n^2}<\frac{b_{n+1}b_{n-1}}{b_n^2}<1+\left(\frac{3(17+12\sqrt{2})}{2(3+2\sqrt{2})^2} + 1 \right) \frac{1}{n^2}, \quad n\geq16,
\]
implying that
\[
1<\frac{b_{n+1}b_{n-1}}{b_n^2}<1+\frac{3}{n^2}, \quad n\geq16.
\]
Let
\[
a_n=\frac{b_n}{n!}, \quad  u_n=\frac{a_{n+1}a_{n-1}}{a_n^2}.
\]
Then we have
\[
\frac{n}{n+1}<u_n<\frac{n}{n+1}\left(1+\frac{3}{n^2}\right), \quad n\geq16.
\]
Let
\[
g_n=\frac{n}{n+1}, \quad \mbox{and} \quad
f_n=\frac{n}{n+1}\left(1+\frac{3}{n^2}\right).
\]
Direct computation shows that
\[
t(g_n,g_{n+1})=\frac{4}{(n+1)(n+2)^2},
\]
which is positive for $n\geq0$.

Similarly, we have
\[
t(g_n,f_{n+1})=\frac{4n^3-21n^2-36n-20}{(n+1)^4(n+2)^2}>0, \quad n>6.
\]

\[
t(f_n,g_{n+1})=\frac{4n^2-21n-9}{n^2(n+1)(n+2)^2}>0, \quad  n>5.
\]

\[
t(f_n,f_{n+1})=\frac{4n^5-24n^4+36n^3+16n^2-48n-144}{n^2(n+1)^4(n+2)^2}>0, \quad  n>7.
\]

Checking the initial values,  we finally get  $\{\frac{b_n}{n!}\}_{n\geq0}$ satisfies the higher order Tur\'an inequality. \qed

\vskip 0.2cm
\noindent{\bf Acknowledgments.}
 We would like to thank the referees for their helpful comments.
  This work was supported by the  National Natural Science Foundation of China (grants 11771330 and 11921001).

\end{document}